\documentclass[12pt]{article}

\title{Ideals of quasi-absorbing elements in semigroups}

\author{
Rico Hager\footnote{formerly Freiberg University of Mining and Technology, Institute for Numerical Mathematics and Optimization, now \href{mailto:ricohager1@gmail.com}{ricohager1@gmail.com}}, Andreas H. Hamel\footnote{Free University of Bolzano-Bozen, Faculty of Economics and Management, \href{mailto:andreas.hamel@unibz.it}{andreas.hamel@unibz.it}}, Frank Heyde\footnote{Freiberg University of Mining and Technology, Institute for Numerical Mathematics and Optimization, \href{mailto:frank.heyde@math.tu-freiberg.de}{frank.heyde@math.tu-freiberg.de}} 
}
\date{{\small \today}}

\usepackage{array} 
\usepackage{verbatim} 

\usepackage{amsmath, amssymb, enumerate, color}
\usepackage{stmaryrd}

\newtheorem{theorem}{Theorem}
\newtheorem{corollary}[theorem]{Corollary}
\newtheorem{remark}[theorem]{Remark}
\newtheorem{lemma}[theorem]{Lemma}
\newtheorem{definition}[theorem]{Definition}
\newtheorem{proposition}[theorem]{Proposition}
\newtheorem{example}[theorem]{Example}

\numberwithin{equation}{section}  
\numberwithin{figure}{section}    
\numberwithin{table}{section}     
\numberwithin{theorem}{section}

\newcommand{\of}[1]{\ensuremath{\left( #1 \right)}}

\newcommand{\cb}[1]{\ensuremath{ \left\{ #1 \right\} }}

\newcommand{\bs}{\backslash}
\newcommand{\vn}{\varnothing}

\newcommand{\pend}{ \hfill $\square$ \medskip}

\renewcommand{\P}{\ensuremath{\mathcal{P}}}

\newcommand{\F}{\ensuremath{\mathcal{F}}}
\newcommand{\G}{\ensuremath{\mathcal{G}}}

\newcommand{\R}{\mathrm{I\negthinspace R}}

\newcommand{\N}{\mathrm{I\negthinspace N}}

\newcommand{\cl}{{\rm cl \,}}

\newcommand{\co}{{\rm conv \,}}

\usepackage[
colorlinks=true, linkcolor=blue, citecolor=blue, urlcolor=blue, filecolor=blue
auskommentieren und unten pdfborder aktivieren]{hyperref}

\definecolor{color0}{gray}{.50}
\definecolor{color1}{rgb}{0,.2,.8}
\definecolor{color2}{rgb}{1,.2,0}
\definecolor{color3}{rgb}{.2,.7,.6}
\newcommand{\f}{\color{color1}}

\begin{document}

\maketitle

\begin{abstract}
Motivated by situations in which the removal of a zero (a.k.a., an absorbing element) from a semigroup yields a subsemigroup with another zero, sets of quasi-zeros (a.k.a., quasi-absorbing elements) are introduced as well as primitive elements, minimal ideals and simple semigroups all with respect to a given ideal. An application to set optimization is discussed.
\end{abstract}



\section{Introduction}

In the past two decades, the theory of set optimization and set-valued (convex) analysis evolved in which the objective functions map into specified subsets of the power set of a linear space \cite{HamelEtAl15Incoll}. These subsets carry the structure of an additive monoid with an absorbing element (a zero in the multiplicative notation in semigroup theory). However, it is a very typical situation that the removal of the (unique) absorbing element makes another element absorbing. The presence of such quasi-absorbing elements has an impact on Green's relations and provides structural information about the monoid. For example, the set of $0$-primitive elements reduces to the quasi-absorbing element (see, e.g., \cite[Theorem 1]{BogdanovicCiricStamenkovic01MM} which becomes trivial if the semigroup includes a quasi-absorbing element). To the best of our knowledge, this situation has not yet been investigated in detail from a semigroup theory point of view. This note aims at filling this gap; in particular, the question is answered which elements have to be eliminated from consideration to end up with non-trivial primitive elements. This leads to concepts such as primitive elements with respect to an ideal and semigroups which are simple with respect to an ideal. It is shown that the set of quasi-absorbing elements is an ideal which in turn provides the set of quasi-primitive elements as well as the new notion of quasi-simple semigroups.

Throughout, $(W, +)$ is a semigroup. The additive notation is used since in set optimization applications, the operation $+$ is a (modification of the) Minkowski addition for subsets of an additive monoid, e.g., the power set of a linear space. 

If $W$ has at least two elements and an element $\alpha \in W$ such that $\alpha + w = w + \alpha = \alpha$ for all $w \in W$, then $\alpha$ is called absorbing (see, e.g., \cite{FlaskaKepka06AUC}; in the multiplicative notation an absorbing element is usually called a zero \cite[p. 2]{Howie95Book}, \cite[p. 3]{Grillet95Book}). An absorbing element, if it exists, is unique. It is called isolated if $(W\bs\{\alpha\}, +)$ is a subsemigroup of $(W, +)$. This is tantamount to saying that there are no two elements $v, w \in W\bs\{\alpha\}$ such that $v + w = \alpha$. 

If $(W, +)$ is a semigroup without (or even with) an absorbing element one can always add (another) one which then becomes isolated. A typical example is $\R\cup\{\pm\infty\}$ with inf-addition, i.e., $(+\infty) + (-\infty) = +\infty$ (all other sums are defined as usual) in which $+\infty$ is isolated absorbing and $-\infty$ becomes isolated absorbing if $+\infty$ is removed; vice versa for the sup-addition featuring $(+\infty) + (-\infty) = -\infty$ and all other sums as before. These two options for adding the infinities are the only ones which keep associativity, i.e., $(+\infty) + (-\infty) = 0$ is not an option. They were introduced by Moreau \cite{Moreau63RCAS} with analytic applications in mind and discussed from an algebraic as well as analytic point of view in \cite{HamelSchrage12JCA}.

An element $e \in W$ is called idempotent if $e + e = e$. The set of idempotents in $W$ is denoted by $E(W)$. The absorbing element is always idempotent as well as the neutral element--if they exist. If $(W, +)$ is commutative, then $(E(W), +)$ is a subsemigroup of $(W, +)$.

\section{Ideals and Green's relations}

In this section, a number of basics concepts and results from Semigroup Theory are collected and shaped in a way suitable for subsequent developments. They can mostly be found in  \cite{Grillet95Book, Howie95Book}.

A set $A \subseteq W$ is called a left (right) ideal in $(W, +)$  if $W + A = \{w + a \mid a \in A, w \in W\} \subseteq A$ ($A + W \subseteq A$). Here and whenever used, the Minkowski operation is extended to empty sets by $\varnothing + A = A + \varnothing = \varnothing$ for all sets $A \subseteq W$.  The case $A = \vn$ is not excluded, i.e., the empty set is always an ideal by definition.

If $A$ is a left as well as a right ideal, it is called a two-sided ideal or just ideal. Clearly, if $A$ is an ideal, then $(A, +)$ is a subsemigroup of $(W, +)$. If $u \in W$ is a fixed element, the set $W_L(u) = (W + \{u\}) \cup \{u\}$ ($W_R(u) = (\{u\} + W) \cup \{u\}$) is the principal left (right) ideal generated by $u$. The principal two-sided ideal of $u \in W$ is the set $W(u) = (W_L(u) + W) \cup W_L(u) = (W + W_R(u)) \cup W_R(u)$.

\begin{lemma}
\label{LemPrincipalIdealRelation}
For all $u \in W$, one has
\begin{multline*}
v \in W_L(u) \, \of{v \in W_R(u), \, v \in W(u)}  \quad \Leftrightarrow \quad \\
	W_L(v) \subseteq W_L(u) \, \of{W_R(v) \subseteq W_R(u), \, W(v) \subseteq W(u)}.
\end{multline*}
\end{lemma}

{\sc Proof.} If $v \in W_L(u)$, then $v = u$ or there is $w \in W$ such that $v = w + u$. Take $v' \in W_L(v)$, i.e., $v' = v \in W_L(u)$ or there is $w' \in W$ such that $v' = w' + v$. In the second case one gets $v' = w' + w + u$, i.e., $v' \in W_L(u)$. Parallel arguments work for $W_R(u)$ and $W(u)$. The reverse implication holds simply because $u \in W_L(u)$ ($u \in W_R(u)$, $u \in W(u)$) for all $u \in W$ by definition. \pend

If $\alpha \in W$ is an absorbing element, then $\alpha \in A$ for any non-empty left or right ideal $A$. In particular, $\alpha \in W_L(u) \cap W_R(u)$ for all $u \in W$ and $W_L(\alpha) = W_R(\alpha) = W(\alpha) = \{\alpha\}$, but there might be another element $\beta \in W\bs\{\alpha\}$ with $\beta \in W_L(u) \cap W_R(u)$ for all $u \in W\bs\{\alpha\}$ which is absorbing on $W\bs\{\alpha\}$--it is precisely this innocent looking fact which motivates the present note.

Green's preorders on a semigroup can be generated by principal ideals. In particular, the $L$-, $R$- and $H$-relation are given by
\begin{align*}
u \leq_L v & \quad \Leftrightarrow \quad W_L(u) \subseteq W_L(v) \; \of{\Leftrightarrow \; u \in W_L(v)}\\
u \leq_R v & \quad \Leftrightarrow \quad W_R(u) \subseteq W_R(v) \; \of{\Leftrightarrow \; u \in W_R(v)} \\
u \leq_H v & \quad \Leftrightarrow \quad W_L(u) \subseteq W_L(v) \; \text{and} \; W_R(u) \subseteq W_R(v) 
	\; \of{\Leftrightarrow \; u \in  W_L(v) \cap W_R(v)}
\end{align*}
for $u, v \in W$. The $L$-, $R$- and $H$-class of $u \in W$ are the sets
\begin{align*}
L(u) & = \cb{v \in W \mid u \leq_L v, \; v \leq_L u} = \cb{v \in W \mid W_L(u) = W_L(v)} \\
R(u) & = \cb{v \in W \mid u \leq_R v, \; v \leq_R u} = \cb{v \in W \mid W_R(u) = W_R(v)} \\
H(u) & = \cb{v \in W \mid u \leq_H v, \; v \leq_H u} = \cb{v \in W \mid W_L(u) = W_L(v), \; W_R(u) = W_R(v)}.
\end{align*}
In a commutative semigroup, one has $L(u) = R(u) = H(u)$ for all $u \in W$.

In the general case, one has the inclusions $H(u) \subseteq L(u) \subseteq W_L(u) \subseteq W(u)$ as well as $H(u) \subseteq R(u) \subseteq  W_R(u) \subseteq W(u)$ for all $u \in W$ which simply follow from the definitions. This raises the question if one can find a set which, when removed from $W(u)$, produces $H(u)$. It turns out that this set always is an ideal.

\begin{lemma}
\label{LemHandASubset}
Let $A \subseteq W$ be a left (right) ideal in $(W, +)$ and $u \in W\bs A$. Then
\[
H(u) \subseteq L(u) \subseteq W_L(u) \bs A \subseteq W(u) \bs A \quad \of{H(u) \subseteq R(u) \subseteq W_R(u) \bs A \subseteq W(u) \bs A}.
\]
\end{lemma}

{\sc Proof.} Assume $A$ is a left ideal and $v \in L(u) \cap A$. Then $W_L(u) = W_L(v)$, hence $u \in W_L(v)$. This gives $u = v \in A$ or $u = w+ v \in A$ for some $w \in W$ which contradicts the assumption in both cases. The proof for right ideals is analogous. The inclusion $W_L(u) \bs A \subseteq W(u) \bs A$ ($W_R(u) \bs A \subseteq W(u) \bs A$) follows from the comment before the lemma. \pend

\begin{theorem}
\label{ThmHandAEq}
Let $A \subseteq W$ be a two-sided ideal in $(W, +)$ and $e \in E(W)\bs A$. Then $H(e) = W(e)\bs A$ if, and only if, $(W(e)\bs A, +)$ is a group.
\end{theorem}

{\sc Proof.} If $H(e) = W(e)\bs A$, then $(W(e)\bs A, +)$ is a group  by Green's theorem. If $(W(e)\bs A, +)$ is a group, then it must be equal to $H(e)$ since $H(e) \subseteq W(e)\bs A$ and $H(e)$ is a maximal subgroup of $(W, +)$ (see \cite[Corollary 1.5, p. 30]{Grillet95Book}). \pend
 
\begin{lemma}(\cite[p. 27]{Grillet95Book})
\label{LemIdempotentIdeals}
Let $e \in E(W)$. 

(1) One has $W_L(e) = W + \{e\}$ and $W_R(e) = \{e\} + W$ as well as
\[
u \in W_L(e) \; (u \in W_R(e)) \quad \Leftrightarrow \quad  u = u + e \; (u = e + u).
\]

(2) One has $W(e) = W_L(e)+W_R(e)$.
\end{lemma}

{\sc Proof.} (1) One gets $W + \{e\} \subseteq W_L(e)$ and $\{e\} + W \subseteq W_R(e)$ from the definitions of $W_L(e)$ and $W_R(e)$. The converse inclusions follow from $e + e = e$.

If $u \in W_L(e) = W + \{e\}$, then there is $w \in W$ such that $u = w+e$ and hence $u + e = w + e + e = w + e = u$. The other direction is immediate. A parallel argument works for $W_R(e)$.

(2) By definition of $W(e)$ and the first part one has
\begin{align*}
 W(e) & = (W_L(e) + W) \cup W_L(e) = (W + \{e\} + W) \cup (W + \{e\})\\
  & = (W + \{e\} + W) \cup (W+\{e\}+\{e\}) = W+\{e\}+W\\
  & = W + \{e\} + \{e\} + W = W_L(e) + W_R(e).
\end{align*}
\pend
 
\begin{theorem}
\label{ThmWnotL}
For all $u \in W$, $W_L(u)\bs L(u)$ is a left and $W_R(u)\bs R(u)$ is a right ideal.
\end{theorem}

{\sc Proof.} If $W_L(u)\bs L(u) = \varnothing$, there is nothing to prove. Otherwise, take $v \in W_L(u) \bs L(u)$. Then $W_L(v) \subseteq W_L(u)$ by Lemma \ref{LemPrincipalIdealRelation}.

It needs to be shown that $W + \{v\} \subseteq W_L(u) \bs L(u)$. One surely has $W + \{v\} \subseteq W_L(v) \subseteq W_L(u)$. If $w + v \in L(u)$ for $w \in W$, then $W_L(u) = W_L(w + v) \subseteq W_L(v)$ by definition of $L(u)$ and Lemma \ref{LemPrincipalIdealRelation}. Altogether, this would imply $W_L(v) = W_L(u)$, i.e., $v \in L(u)$ contradicting the assumption.

The argument for $W_R(e)\bs R(e)$ runs parallel. \pend

\begin{corollary}
\label{CorWnotH}
If $(W, +)$ is a commutative semigroup and $u \in W$, then $W(u)\bs H(u)$ is an ideal in $(W, +)$.
\end{corollary}

\begin{example}
\label{ExNonCommutative}
The commutativity is indeed necessary as can be seen from the following simple example.
\begin{center}
\begin{tabular}{ c | c c c c}
 $+$ & $e$ & $f$\\ 
\hline
 $e$ & $e$ & $e$\\ 
 $f$ & $f$ & $f$\\
\end{tabular}
\end{center}
For $x \in \{e, f\} = W$ one gets $W_L(x) = W(x) = W$, $W_R(x) =\{x\} = H(x)$, and $W(e)\bs H(e) = \{f\}$ is not an ideal.
\end{example}

Theorem \ref{ThmWnotL} and Corollary \ref{CorWnotH} justify that $A$ in the previous results is assumed to be an ideal.

\section{Quasi-absorbing and quasi-primitive elements}

It is well-known that $\leq_H$ is a partial order on $E(W)$ with $e \leq_H f$ if, and only if, $e + f = f + e = e$ for $e, f \in E(W)$ (compare \cite[Sec. 3.2]{Howie95Book}). This partial order is sometimes called the Rees order (see \cite[Sec. II.1, p. 27]{Grillet95Book}). The set of idempotents dominated by $e \in E(W)$ with respect to $\leq_H$ is denoted by $E^\leq(e) = \{f \in E(W) \mid f \leq_H e\}$. For $f, g \in E^\leq(e)$ one has $f \leq_H g$ if, and only if, $f + g = g + f = f$ if, and only if, $E^\leq(f) \subseteq E^\leq(g)$. 

The following definition formalizes the intuitive idea to successively remove absorbing elements one by one. It starts with the absorbing element $\alpha$ in $(W, +)$, then the one in $(W\backslash\{\alpha\}, +)$ and so forth until no more absorbing elements can be found. The blueprint example is $\R\cup\{\pm\infty\}$ with inf-addition: one would remove the absorbing element $+\infty$ from $\R\cup\{\pm\infty\}$ and then $-\infty$ from $\R\cup\{-\infty\}$. However, in some cases no elements can be removed in this way although there are candidates for "quasi-absorbing" element.

\begin{definition}
\label{DefAs(W)}
The set of stepwise quasi-absorbing elements of a semigroup $(W,+)$ is defined by
\[
A_s(W) = \{\alpha\in W \mid \exists n\in \mathbb{N} \colon \alpha \in A_n\},
\]
where $A_0 = \emptyset$, $\bar{A}_n = \bigcup_{i=1}^{n}A_i$ and
\[
A_n = \{\alpha \in W \bs \bar{A}_{n-1} \mid \forall w\in W\bs \bar{A}_{n-1} \colon w + \alpha = \alpha + w = \alpha\}
\]
for $n \geq 1$.
\end{definition}

\begin{example}
\label{ExNatNumbersMax} The commutative semigroup of the natural numbers $W= \mathbb N$ with the max-operation for two numbers is an example for $A_s(W) = \varnothing$ since there is no absorbing element. On the other hand, one has $\max\{m, n\} = n$ for all $m \leq n$, $m, n \in \N$, which could be understood as: every number absorbs all numbers less than or equal to it. One may ask if there is a concept of quasi-absorbing elements which would produce that all numbers are quasi-absorbing in this case.
\end{example}

The elements of $A_s(W)$ are idempotent. Indeed, the sets $A_n$ are singletons or empty: if $\alpha_1, \alpha_2 \in A_n$, then $\alpha_1 + \alpha_2 = \alpha_2 + \alpha_1= \alpha_2$ as well as $\alpha_2 + \alpha_1 = \alpha_1 + \alpha_2= \alpha_1$ which gives $\alpha_1 = \alpha_2$. Consequently, $\alpha + \alpha = \alpha$ for $A_n = \{\alpha\}$. Clearly $A_n = \varnothing$ implies $A_m = \varnothing$ for $m >n$. Let $\alpha_n$ denote the single element in $A_n$ whenever $A_n \neq \varnothing$. Then one has $\alpha_{n+1} + \alpha_n = \alpha_n + \alpha_{n+1} = \alpha_n$ which means $\alpha_n \leq_H \alpha_{n+1}$ for $n = 1,2, \ldots$ whenever $A_{n+1} \neq \varnothing$. Therefore, $\leq_H$ is a total order on $A_s(W)$ where the terminology of \cite[Sec. 1.3]{Howie95Book} is used. This motivates the following concepts.

With the notation
\[
E^{lin}(W)  = \cb{\alpha \in E(W) \mid \; \leq_H \;  \text{is a total order on} \; E^\leq(\alpha)}
\]
one always has $\alpha \in E^{lin}(W)$ for an absorbing element $\alpha \in W$ since $E^\leq(\alpha) = \{\alpha\}$. Define the set
\[
A_{fin}(W) = \{\alpha \in E^{lin}(W) \mid E^\leq(\alpha) \; \text{is finite}, \; \forall w \in W\bs E^\leq(\alpha) \colon w + \alpha = \alpha + w = \alpha\}.
\]

\begin{lemma}
\label{LemAWChar}
For $\beta\in W$ one has
\[
\beta \in A_{fin}(W) \quad \Leftrightarrow \quad \exists \alpha \in A_{fin}(W) \colon \beta \in E^\leq(\alpha).
\]
\end{lemma}

{\sc Proof.} If $\beta \in A_{fin}(W)$, then one can take $\alpha = \beta$ on the right hand side of the implication since $\beta \in E^\leq(\beta)$. Conversely, assume $\beta \in E^\leq(\alpha)$ for some $\alpha \in A_{fin}(W)$. Then $\beta$ is idempotent since $E^\leq(\alpha) \subseteq E(W)$. Moreover, $\leq_H$ is a total order on $E^\leq(\beta)$ since  $E^\leq(\beta) \subseteq E^\leq(\alpha)$ and it is one on $E^\leq(\alpha)$.

Take $w \in W \bs E^\leq(\beta)$. It needs to be verified that $w+ \beta = \beta + w = \beta$. If $w \in W \bs E^\leq(\alpha)$ then
\[
\beta + w = (\beta + \alpha) + w = \beta + (\alpha + w) = \beta + \alpha = \beta
\]
and
\[
w + \beta = w + (\alpha + \beta) = (w + \alpha) + \beta = \alpha + \beta = \beta.
\]
Since $\leq_H$ is a total order on $E^\leq(\alpha)$, $\beta \leq_H w$ or $w \leq_H \beta$ holds for all $w, \beta \in E^\leq(\alpha)$. If $w \in E^\leq(\alpha)\bs E^\leq(\beta)$, then the second case is excluded, hence $\beta \leq_H w$ which yields the claim.
\pend

\begin{lemma}
\label{LemAWProperties}
(1) If $\alpha, \beta \in A_{fin}(W)$,  then $\alpha + \beta = \beta + \alpha \in \{\alpha, \beta\}$.

(2) $A_{fin}(W)$ is an ideal in $(W, +)$.

(3) $A_{fin}(W) \subseteq W_L(u) \cap W_R(u)$ for all $u \in W\bs A_{fin}(W)$. 
\end{lemma}

{\sc Proof.} (1) This is clearly true for $\alpha = \beta$. Assume $\alpha \neq \beta$. Then either $\beta \in W \bs E^\leq(\alpha)$ and hence $\alpha + \beta =  \beta + \alpha = \alpha$ by definition of $A_{fin}(W)$, or $\beta \in E^\leq(\alpha)$ and hence $\beta + \alpha = \alpha + \beta = \beta$. 

(2) If $A_{fin}(W) = \vn$, then there is nothing to prove. Otherwise, take $\alpha \in A_{fin}(W)$ and $w \in W$. If $w \in E^\leq(\alpha)$, then Lemma \ref{LemAWChar} yields $w \in A_{fin}(W)$ and (1) gives $w + \alpha = \alpha + w \in A_{fin}(W)$. If $w \in W\bs E^\leq(\alpha)$, then $w + \alpha = \alpha + w = \alpha \in A_{fin}(W)$ by definition of $A_{fin}(W)$.

(3) If $A_{fin}(W) = \vn$, the statements holds trivially. If $A_{fin}(W) \neq \vn$, one has $u + \alpha = \alpha + u = \alpha$ and hence $\alpha \in W_L(u) \cap W_R(u)$ for all $u \in W \bs A_{fin}(W)$ and all $\alpha \in A_{fin}(W)$.  
\pend

\begin{theorem}
\label{ThmAs=Alin}
One has $A_s(W)=A_{fin}(W)$.
\end{theorem}

{\sc Proof.} Let $\alpha \in A_s(W)$. Then, there is $n \in \mathbb{N}$ such that $\{\alpha\} = A_n$ since the sets $A_n$ are singletons. Therefore the abbreviation $\alpha_n = \alpha$ is used in the following. One also has $\bar{A}_n = \{\alpha_1, \ldots, \alpha_n\}$. If $1 \leq m < n$, the definition of $\alpha_m$ gives $\alpha_n + \alpha_m = \alpha_m + \alpha_n = \alpha_m$, i.e., $\alpha_m \leq_H \alpha_n$, and this is also true for $n = m$ since $\alpha_n$ is idempotent.

Moreover, there cannot be another idempotent $\beta \in E(W) \bs \{\alpha_1, \ldots, \alpha_n\}$ such that $\alpha_n + \beta = \beta + \alpha_n = \beta$ since this would contradict $\alpha_n \in A_n$. Therefore, $E^\leq(\alpha_n)=\{\alpha_1, \ldots, \alpha_n\} = \bar A_n$ which is a finite set. This implies  $w + \alpha_n = \alpha_n + w = \alpha_n$ for all $w \in W\bs E^\leq(\alpha_n)$. It remains to show that $\leq_H$ is a total order on $E^\leq(\alpha_n)$ and thus $\alpha_n\in E^{lin}(W)$. Because of $E^\leq(\alpha_n) \subseteq E(W)$,  $\leq_H$ is a partial order on $E^\leq(\alpha_n)$. If one takes $\alpha_i, \alpha_j \in \{\alpha_1, \ldots, \alpha_n\} = E^\leq(\alpha_n)$ with $1 \leq i \leq j \leq n$, then $a_i \leq_H \alpha_j$ holds true which means that $\leq_H$ must be total. Altogether, $\alpha = \alpha_n\in A_{fin}(W)$.

To show the converse inclusion, take $\alpha \in A_{fin}(W)$. The set $E^\leq(\alpha)$ is finite and totally ordered by $\leq_H$, i.e., one has $\alpha_1 \leq_H \ldots \leq_H \alpha_n$ for its $n$ elements where $\alpha_n = \alpha$ and $\alpha_i \leq_H \alpha_j$ if, and only if,  $i \leq j$ for $i, j \in \{1, \ldots, n\}$. This implies $E^\leq(\alpha_i) \subseteq E^\leq(\alpha_j)$ if, and only if, $i \leq j$  for $i, j \in \{1, \ldots, n\}$. Therefore $E^\leq(\alpha_j) = \{\alpha_1, \ldots, \alpha_j\}$ for all $j \in \{1, \ldots, n\}$. Lemma \ref{LemAWChar} yields $\alpha_j \in A_{fin}(W)$ for all $j \in \{1, \ldots, n\}$. 

By induction it will be shown that $A_j = \{\alpha_j\}$ with $A_j$ from Definition \ref{DefAs(W)} for $j \in \{1, \ldots, n\}$. For $\alpha_1 \in A_{fin}(W)\subseteq E(W)$ one has $\alpha_1 \in A_1$ since $E^\leq(\alpha_1) = \{\alpha_1\}$ and thus $\alpha_1 \in A_{fin}(W)$ means that $\alpha_1$ is absorbing in $(W, +)$. Since $A_1$ is a singleton, $A_1 = \{\alpha_1\}$.

Assume $1< j \leq n$ and $A_i = \{\alpha_i\}$ for $i \in \{1, \ldots, j-1\}$. Then $\alpha_j \in W\bs \bar{A}_{j-1} = W\bs \{\alpha_1, \ldots, \alpha_{j-1}\}$ and $\alpha_j \in A_{fin}(W)$ implies $w + \alpha_j = \alpha_j + w = \alpha_j$ for all $w \in W\bs E^\leq(\alpha_j) = W\bs\{\alpha_1, \ldots, \alpha_j\}$. This equation is also true for $w = \alpha_j$ since $\alpha_j \in E(W)$, hence for all $w \in W\bs\{\alpha_1, \ldots, \alpha_{j-1}\} = W\bs \bar{A}_{j-1}$. Therefore $A_j = \{\alpha_j\}$ and consequently, $\{\alpha\} = \{\alpha_n\} = A_n$ which yields $\alpha \in A_s(W)$.
\pend

The alternative representation of $A_s(W)$ in Theorem \ref{ThmAs=Alin} is generalized in the next definition. The requirement that the sets $E^\leq(\alpha)$ are finite is dropped.

\begin{definition}
The elements of the set
\[
A(W) = \cb{\alpha \in E^{lin}(W) \mid \forall w \in W \bs E^\leq(\alpha) \colon w + \alpha = \alpha + w = \alpha}
\]
are called quasi-absorbing.
\end{definition}

If $\alpha \in W$ is absorbing, then $\alpha \in A(W)$ with $E^\leq(\alpha) = \{\alpha\}$. On the other hand, $A(W) \neq \vn$ does not imply that $(W, +)$ has an absorbing element: consider the semigroup $(\R_+\bs\{0\}, \wedge)$ with $A(W) = \R_+\bs\{0\}$.
 
\begin{remark}
\label{RemAWProperties}
An inspection of the proofs of Lemma \ref{LemAWChar} and Lemma \ref{LemAWProperties} shows that they remain valid for $A(W)$ since the finiteness assumption for the sets $E^\leq(\alpha)$ is not used in those proofs.
\end{remark}

\begin{remark}
\label{RemAWInclusion}
The inclusion $A_s(W) = A_{fin}(W) \subseteq A(W)$ follows from the definitions. Moreover, if $W$ is finite, then $A(W) = A_{fin}(W)$. However, the strict inclusion might occur: in Example \ref{ExNatNumbersMax}, one has $A_s(W) = \varnothing$ and $A(W)= \mathbb N$ since $E^\leq(n) = \{m \in \mathbb N \mid \max\{m, n\} = m\} = \{m \in \mathbb N \mid m \geq n\}$. 

One can even have $\varnothing \neq A_{fin}(W) \subsetneq A(W)$: let $W=\{1-\frac{1}{n} \mid n \in \mathbb N\backslash\{0\}\} \cup\{1+\frac{1}{n} \mid \in \mathbb N\backslash\{0\}\}$ which is a semigroup with the max operation. Then, one has $A(W) = W$ and $A_{fin}(W)=\{1+\frac{1}{n} \mid n \in\mathbb N\backslash\{0\}\}$. This example also shows that there might not be an $\alpha \in A(W)$ such that $E^\leq(\alpha)=A_{fin}(W)$.
\end{remark}
 
One may also note that (1) of Lemma \ref{LemAWProperties} means that $(A_s(W)= A_{fin}(W), +)$ and $(A(W), +)$ are commutative subsemigroups of $(W, +)$ with only idempotent elements and thus can be identified with lower semilattices. 

\begin{corollary}
Let $e \in E(W)\bs A(W)$ ($e \in E(W)\bs A_s(W)$). Then  $H(e) = W(e)\bs A(W)$ ($H(e) = W(e)\bs A_s(W)$) holds true if, and only if, $(W(e)\bs A(W), +)$ ($(W(e)\bs A_s(W), +)$) is a group with neutral element $e$.
\end{corollary}

{\sc Proof.} This is a direct consequence of Theorem \ref{ThmHandAEq} and Lemma \ref{LemAWProperties}, (2) in connection with Remark \ref{RemAWProperties}. \pend

Primitive elements are minimal idempotents with respect to the Rees order within the set of non-absorbing idempotents (see \cite[Section 3.2]{Howie95Book}).  They are relevant for the definition of completely 0-simple semigroups and thus for the Rees Theorem \cite[Theorem 3.2.3]{Howie95Book}. However, after the removal of an absorbing idempotent another one might become absorbing and thus the only primitive element (and also the unique minimum) in the remaining set.  This might render some results trivial and thus motivates the following definition.

\begin{definition}
Let $A \subseteq W$ be an ideal in $(W, +)$. An element $e \in E(W) \bs A$ is called $A$-primitive if it satisfies
\[
f \in E(W), \, f \leq_H e \quad \Rightarrow \quad f = e \; \text{or} \; f \in A.
\]
The set of $A$-primitive elements is denoted by $P(A)$. An $A(W)$-primitive idempotent is called quasi-primitive and $P(A(W))$ is the set of quasi-primitive elements.
\end{definition}

If $A, B \subseteq W$ are two ideals with $A \subseteq B$ and $e \in E(W) \bs A$ is $A$-primitive, then there are two cases possible: first, $e \not\in B$ in which case $e$ is also $B$-primitive since $A \subseteq B$; secondly, $e \in B \bs A$. 

\begin{theorem}
\label{ThmQAPrimitive}
Let $A \subseteq W$ be an ideal in $(W, +)$ and $e \in P(A)$. Then
\[
H(e) = \of{W_L(e) \cap W_R(e)} \bs \bigcup_{f \in P(A)}\of{W_L(f) \cap W_R(f)}\bs H(f).
\]
\end{theorem}

{\sc Proof.} One always has the inclusion
\begin{align*}
H(e) & = W_L(e) \cap W_R(e) \cap H(e) \\
	& = \of{W_L(e) \cap W_R(e)} \bs \of{W_L(e) \cap W_R(e) \bs H(e)} \\
	& \supseteq \of{W_L(e) \cap W_R(e)} \bs \bigcup_{f \in P(A)}\of{W_L(f) \cap W_R(f)}\bs H(f).
\end{align*}
To prove the reverse inclusion, take $w \in H(e)$. Then $w \in W_L(w) = W_L(e)$ and $w \in W_R(w) = W_R(e)$. Assume $w \in \bigcup_{f \in P(A)}\of{W_L(f) \cap W_R(f)}\bs H(f)$, i.e., there is $g \in P(A)$ with $w \in \of{W_L(g) \cap W_R(g)}\bs H(g)$. Clearly $e \neq g$, $w \neq g$ since $w \in H(e)$, $w \not\in H(g)$. Since $w \in W_L(g)$, one has $W_L(e) = W_L(w) \subseteq W_L(g)$ by Lemma \ref{LemPrincipalIdealRelation}. A parallel argument ensures $W_R(e) \subseteq W_R(g)$. Altogether, this implies $e \leq_H g$, $e \neq g$ which contradicts $g \in P(A)$. This proves the remaining inclusion. \pend

\begin{corollary}
\label{CorAPrimitive}
Let $(W, +)$ be a commutative semigroup and $A \subseteq W$ an ideal in $(W, +)$. If $e \in P(A)$, then 
\[
H(e) = W(e) \bs \bigcup_{f \in P(A)} W(f)\bs H(f).
\]
\end{corollary}

Example \ref{ExNonCommutative} also shows that commutativity is necessary: in this case, one has  $H(x) \neq W(x) \bs W = \vn$ for $x \in W$.

\begin{corollary}
\label{CorAPrimitiveIdeal}
Let $(W, +)$ be a commutative semigroup and $A \subseteq W$ an ideal in $(W, +)$. Then, the set
\[ 
\bigcup_{f \in P(A)} W(f)\bs H(f).
\]
is an ideal in $(W, +)$.
\end{corollary}

The benefit of $A(W)$ becomes apparent when considering primitive elements. The following theorem shows that the set of $A_{fin} (W)$-primitive elements will be in some sense trivial if $A_{fin} (W)$ is different from $A(W)$.

\begin{theorem}
\label{ThmFinQuasiPrimitives}
If $A_{fin}(W) \subsetneq A(W)$ then $P(A_{fin}(W)) = \varnothing$ or $P(A_{fin}(W)) = \{e\}$ with $e\in A(W)$.
\end{theorem}

{\sc Proof.} If $A_{fin}(W) \subsetneq A(W)$ then there is some $\alpha \in A(W)\bs A_{fin}(W)$. The set $E^\le(\alpha)\bs A_{fin}(W)$ is nonempty and totally ordered with respect to $\leq_H$. If it has a least element $e$, then $E^\leq(e)\bs A_{fin}(W)=\{e\}$, hence $f \in E(W)$, $f \leq_H e$ implies $f=e$ or $f \in A_{fin}(W)$. This means $e \in P(A_{fin}(W))$. Moreover, $e\in E^\leq(\alpha)$ implies $e \in A(W)$ due to Lemma \ref{LemAWChar} and Remark \ref{RemAWProperties}. Assume, that there would be another element $f \in P(A_{fin}(W))$. Then, $f \not\in A_{fin}(W)$ and $e \not\leq_H f$ must be true. But, $e \in A(W)$ implies $e \leq_H f$ whenever $f\ in E(W)\bs E^\leq(e)$, hence $e\not\leq_H f$ is only possible for $f \in E^\leq(e)\bs \{e\} \subseteq A_{fin}(W)$ which would be a contradiction. If $E^\leq(\alpha)\bs A_{fin}(W)$ has no least element, then for all $e \in E^\leq(\alpha)\bs A_{fin}(W)$ there is some $f \in E^\leq(\alpha)\bs A_{fin}(W)$ with $f\leq_H e$ and $f \neq e$. Moreover, for all $e \in E(W)\bs E^\leq(\alpha)$ one has $\alpha \leq_H e$. Hence $P(A_{fin}(W)) = \varnothing$. 
\pend

\begin{example}
\label{ExPFinNotPAW}
Consider the semigroup $(W,*)$ with $W=D\cup\{e,f\}$, $D=\{\frac{1}{n} \mid n \in \N \bs\{0\}\}$ where the binary operation $*$ is defined as $d_1*d_2=\min\{d_1, d_2\}$ for $d_1, d_2 \in D$ and according to the following table for $e,f$ and elements $d\in D$.
\begin{center}
\begin{tabular}{ c | c c c }
$*$ & $d$ & $e$ & $f$\\
\hline
$d$ & $d$ & $d$ & $d$\rule{0pt}{3ex}\\
$e$ & $d$ & $e$ & $1$\rule{0pt}{3ex}\\
$f$ & $d$ & $1$ & $f$\rule{0pt}{3ex}\\ 
\end{tabular}
\end{center}
Then $A_{fin}(W) = \varnothing$, $P(A_{fin}(W)) = \varnothing$ and $A(W)=D$, $P(A(W)) = \{e,f\}$.
\end{example}

In Example \ref{ExPFinNotPAW}, $P(A_{fin} (W))$ is trivial while $P(A(W))$ is not. Of course it may happen that both sets are trivial, even if $A_{fin} (W)$ and $A(W)$ do not coincide. A particular case is $A(W)=W$ and hence $P(A(W)) = \varnothing$ as in the following example.

\begin{example}
\label{ExTrivialWAs}
Define $W = D_1 \cup D_2 \subset \N$ with $D_1 = \cb{1 - \frac{1}{n} \mid n \in \N\bs\{0\}}$ and $D_2 = \cb{2 - \frac{1}{n} \mid n \in \N\bs\{0\}}$. Considering $(W, \min)$ one gets $A_{fin}(W) = D_1$ with $P(A_{fin}(W)) = \{1\}$ as well as $A(W) = W$ with $P(A(W)) = \varnothing$. This illustrates the second case in Theorem \ref{ThmFinQuasiPrimitives}. Note that $\alpha = 0$ is the absorbing element in $(W, \min)$.
\end{example}

\section{Simple semigroups with respect to ideals}

The next two definitions give generalizations of $0$-minimal ideals and $0$-simple semigroups, compare \cite[p. 44f]{Grillet95Book}, \cite[p. 66f]{Howie95Book}.

\begin{definition} 
\label{DefAMinimalIdeals}
Let $A \subseteq W$ be an ideal of $(W,+)$. An ideal $D$ in $(W,+)$ is called $A$-minimal if $D \not\subseteq A$ and $D$ does not properly contain an ideal in $(W,+)$ other than subsets of $A$. An $A(W)$-minimal ideal is called quasi-minimal.
\end{definition}

\begin{definition} 
\label{DefASimpleSemigroups}
Let $A \subsetneq W$ be an ideal in $(W,+)$. The semigroup $(W,+)$ is called $A$-simple if the only ideals in $(W,+)$ are $W$ and subsets of $A$. An $A(W)$-simple semigroup is called quasi-simple.
\end{definition}

If $\alpha \in W$ is an absorbing element, then the choice $A = \{\alpha\}$ produces $\alpha$-minimal ideals and $\alpha$-simple semigroups usually called $0$-minimal and $0$-simple in the literature.

\begin{remark}
\label{RemSubsemigroup} 
Let $A, D \subseteq W$ be ideals in $(W, +)$ with $A \cap D \subsetneq D$. Then $D$ is an $A$-minimal ideal in $(W,+)$ if $(D,+)$ is an $(A \cap D)$-simple subsemigroup of $(W,+)$. Indeed, the condition $A \cap D \subsetneq D$ is tantamount to $D \not\subseteq A$ and an ideal $E$ in $(W, +)$ with $E \subseteq D$ is also an ideal in $(D, +)$.
\end{remark}

\begin{theorem}
\label{ThmAMinimalIdealsChar}
Let $(W, +)$ be a commutative semigroup, $A \subseteq W$ an ideal in $(W, +)$ and $e \in E(W)\bs A$. Then, the following statements are equivalent:

(1) $H(e)=W(e)\bs A$;

(2) $W(e) \bs A$ is a group;

(3) $W(e)$ is an $A$-minimal ideal of $W$;

(4) $W(e)$ is an $A \cap W(e)$-simple semigroup. 
\end{theorem}

\begin{remark}
\label{RemENotInA}
Note that $e \not\in A$ implies $e \in W(e)\bs A$, hence $A\cap W(e) \subsetneq W(e)$.
\end{remark}

{\sc Proof.} (1)$\Rightarrow$(2): This implication directly follows from \cite[Proposition II.1.4]{Grillet95Book} (which can be seen as a version of Green's Theorem).

(2)$\Rightarrow$(4): Suppose there is an ideal $D$ of $W(e)$ such that $D\not\subseteq A$ and $D\neq W(e)$. Take $u \in D \bs A\subseteq W(e) \bs A$. Since $e \in W(e)\bs A$ and $W(e) \bs A$ is a group, there is $v \in W(e) \bs A$ such that $u + v = e$. Hence $e \in D$, since $D$ is an ideal of $W(e)$. 

Let $w\in W(e)$ be arbitrary. By Lemma \ref{LemIdempotentIdeals} (2)  there are $u,v\in W$ such that
$w=u+v$ with $u\in W_L(e)\subseteq W(e)$, $v \in W_R(e) \subseteq W(e)$. One gets $u=u + e \in W(e)+D\subseteq D$ from Lemma \ref{LemIdempotentIdeals} (1). Hence $w \in D + W(e) \subseteq D$. This shows $W(e)\subseteq D$, hence $D = W(e)$ in contradiction to $D\neq W(e)$.

(4)$\Rightarrow$(3): This implication is immediate from Remark \ref{RemSubsemigroup}.

(3)$\Rightarrow$(1): By Lemma \ref{LemHandASubset}, one has $H(e) \subseteq W(e) \bs A$. It remains to show that $W(e) \bs A \subseteq H(e)$. If $u \in W(e) \bs A$, then $W(u) \subseteq W(e)$ since $W(e)$ is an ideal of $W$. Moreover, $W(u) \not\subseteq A$, since $u \in W(u)$, but $u \not\in A$. Now, the $A$-minimality of $W(e)$ implies $W(u)=W(e)$ and, consequently, $u \in H(e)$. 
\pend

\begin{remark}
(1)$\Leftrightarrow$(2) is true also for non-commutative semigroups (see Theorem \ref{ThmHandAEq}) as well as the implications (4)$\Rightarrow$(3) and (2)$\Rightarrow$(4). (3)$\Rightarrow$(1) becomes a problem since in Example \ref{ExNonCommutative} one has $W(e) = W(f) = W$, but $f \not\in H(e) = \{e\}$ and thus the last conclusion in the proof cannot be drawn. Moreover, Example \ref{ExNonCommutative} also shows that this implication is not true in general: in the example, $W(e)$ is an $\vn$-minimal ideal with $H(e) = \{e\} \neq \{e, f\} = W(e)\bs \vn$. It is also possible to add an isolated absorbing element $\alpha$ to $W= \{e,f\}$ and get the same effect for the ideal $A = \{\alpha\}$.
\end{remark}

\begin{theorem}
\label{ThmARepGroup}
Let $(W,+)$ be a commutative semigroup, $A$ be an ideal  in $(W, +)$ and $e\in P(A)$. Then
\[
(W(e) \bs A, +) \;\text{is a group} \quad\Leftrightarrow\quad W(e) \cap A = W(e) \cap \bigcup_{f\in P(A)} W(f) \bs H(f).
\]
\end{theorem}

{\sc Proof.}  
By Corollary \ref{CorAPrimitive} we have 
\[
H(e) = W(e) \bs \bigcup_{f \in P(A)} W(f)\bs H(f).
\]
Since $e\in P(A)$ implies $e\in E(W)\bs A$, one obtains
\begin{align*}
(W(e)\bs A, +) \; \text{is a group} \quad & \Leftrightarrow\quad W(e)\bs A = H(e) = W(e) \bs \bigcup_{f \in P(A)} W(f)\bs H(f) \\
                                             & \Leftrightarrow\quad W(e) \cap A = W(e) \cap \bigcup_{f \in P(A)} W(f) \bs H(f)
\end{align*}
with the help of Theorem \ref{ThmHandAEq}.
\pend

\begin{corollary}
Let the assumptions of Theorem \ref{ThmARepGroup} be satisfied and additionally $A\subseteq W(e)$ and $\bigcup_{f\in P(A)} W(f) \bs H(f)\subseteq W(e)$. Then
\[
(W(e) \bs A, +) \;\text{is a group} \quad\Leftrightarrow\quad A = \bigcup_{f\in P(A)} W(f) \bs H(f).
\]
\end{corollary}

One may ask for a sufficient condition for the last assumption in the previous corollary. It can be formulated as a condition to an ideal $A \subseteq W$.

\begin{definition}
\label{DefBottleneckIdeal} An ideal $A \subseteq W$  in $(W, +)$ is called a bottleneck ideal if for all ideals $D \subseteq W$ in $(W, +)$ one has $D \subseteq A$ or $A \subseteq D$ (or both).
\end{definition}

\begin{lemma}
If $A$ is a bottleneck ideal in $(W, +)$ and $e \in E(W)$, then $e \not\in A$ implies $A \subsetneq W(e)$. 
\end{lemma}

{\sc Proof.} If $e \not\in A$, then $W(e) \not\subseteq A$, hence $A \subseteq W(e)$ since $A$ is a bottleneck ideal. Moreover, $A$ is a proper subset of $W(e)$ since $e\in W(e)\bs A$.
\pend

\begin{proposition}
\label{PropAWBottle}
The sets $A(W)$ and $A_s(W)$ are bottleneck ideals in $(W, +)$.
\end{proposition}

{\sc Proof.} Assume $D \not\subseteq A(W)$, i.e., there is $u \in D$ with $u \not\in A(W)$. Then $u \not\in E^{\le}(a)$ for all $a \in A(W)$ by Lemma \ref{LemAWChar} and hence $a = a + u = u + a$ according to the definition of $A(W)$, i.e., $a \in D$ for all $a \in A(W)$. The proof for $A_s(W)$ runs in the same way.
\pend

\begin{example}
\label{ExBottleneckIdeal1}
The set $A = \{w, \alpha\}$ is a bottleneck ideal in the semigroup given by the following addition table:
\begin{center}
\begin{tabular}{ c | c c c c}
 $+$ & $u$ & $v$& $w$ & $\alpha$\\ 
\hline
 $u$ & $u$ & $v$ & $w$ & $\alpha$\\ 
 $v$ & $v$ & $u$ & $w$ & $\alpha$\\
 $w$ & $w$ & $w$ & $\alpha$ & $\alpha$\\
 $\alpha$ & $\alpha$ & $\alpha$ & $\alpha$ & $\alpha$\\
\end{tabular}
\end{center}
One also has $A(W) = \{\alpha\}$. Thus, a bottleneck ideal is not necessarily a subset of $A(W)$. On the other hand, $E^\leq(a)$ is a bottleneck ideal for each $a \in A(W)$ in every semigroup $(W, +)$.
\end{example}

\begin{example}
\label{ExBottleneckIdeal2}
The commutative semigroup $(\R_+, +)$ of non-negative numbers with the usual addition only has bottleneck ideals. Indeed, an ideal is a set $A \subseteq \R_+$ satisfying $A + \R_+ \subseteq A$, i.e., $A$ is of the form $(a, \infty)$, $[a, \infty)$ with $a \in \R_+$. The collection of these sets is totally ordered by inclusion. The only idempotent is 0 and $A(W) = \vn$.
\end{example}

\begin{remark}
In \cite[Definition 1]{Abrhan94MS} (see also \cite{Abrhan97MB}), an ideal $D \subseteq W$ is called minimal with respect to a subset $\vn \neq B \subseteq W$ if $D \cap B \neq \vn$ and there is no ideal $D' \subsetneq D$ with $D' \cap B \neq \vn$. Thus, an ideal $D \subseteq W$ is minimal with respect to a subset $\vn \neq W\bs A$, $A \neq W$, if $D \cap (W \bs A) \neq \vn$ and there is no ideal $D' \subsetneq D$ with $D' \cap W \bs A \neq \vn$. The last condition is equivalent to $D' \not\subseteq A$.

Thus, every $A$-minimal ideal is also minimal with respect to $W\bs A$ in the sense of Abrhan \cite[Definition 1]{Abrhan94MS}. The converse is not generally true. However, if $A$ is a bottleneck ideal, then an ideal $D$ is $A$-minimal if, and only if, it is minimal with respect to $W\bs A$ in the sense of Abrhan \cite[Definition 1]{Abrhan94MS}. In particular, an ideal is $A(W)$-minimal if, and only if, it is minimal with respect to $W\bs A(W)$ in the sense of \cite{Abrhan94MS}.
\end{remark}

\section{Applications in set optimization}

The following construction furnishes semigroups with relevance in set-valued convex analysis and for optimization problems with a set-valued objective as surveyed in \cite{HamelEtAl15Incoll}. This approach has applications in Mathematical Finance \cite{HamelHeyde10SIFIN}, Statistics \cite{HamelKostner18JMVA}, Game Theory \cite{HamelLoehne18MMOR} and Vector Optimization \cite{Loehne11Book}, for example.

Let $d \geq 1$ be a natural number. A convex cone $D \subseteq \R^d$ is a set which satisfies $s z \in D$ whenever $s > 0$, $z \in D$ and $z_1+z_2 \in D$ whenever $z_1, z_2 \in D$. Let $C \subseteq \R^d$ be a (fixed) closed convex cone and 
\[
\G(\R^d, C) = \cb{A \subseteq \R^d \mid A = \cl\co(A + C)}
\]
the collections of subsets of $\R^d$ which are closed, convex and invariant under addition of $C$. Here, $A + B = \{a + b \mid a, \in A, b \in B\}$ is the elementwise addition of sets with $\vn + A = A + \vn = \vn$ for each $A \subseteq \R^d$. Moreover, $\cl$ and $\co$ denote the topological closure and the convex hull, respectively, of a set in $\R^d$. The addition in $\G(\R^d, C)$ is given by $A \oplus B = \cl(A+B)$ which turns $(\G(\R^d, C), \oplus)$ into a commutative monoid with neutral element $C$ and (isolated) absorbing element $\vn$. Note that $(\G(\R^d, C)\bs\{\vn\}, \oplus)$ is a commutative monoid as well with absorbing element $\R^d$ which is not isolated in general.

The dual cone of $C$ is the set $C^+ = \{y \in \R^d \mid \forall z \in C \colon y^\top z \geq 0\}$ which also is a closed convex cone. If $y \in \R^d\bs\{0\}$, $h^+_y = \{z \in \R^d \mid y^\top z \geq 0\}$ is the closed (homogeneous) halfspace with normal $y$, and one has $h^+_y \in \G(\R^d, C)$ if, and only if, $y \in C^+$.

\begin{lemma}
\label{LemElementaryInG}
(1) A set $A \in \G(\R^d, C)$ is idempotent in $(\G(\R^d, C), \oplus)$ if, and only if, it is a closed convex cone with $A \supseteq C$ or $\vn$.

(2) On $E(\G(\R^d, C))\bs\{\vn\}$, the Rees order $\leq_H$ coincides with $\supseteq$.

(3) $\vn, \R^d$ are quasi-absorbing in $(\G(\R^d, C), \oplus)$.

(4) Half-spaces $h^+_y$ with $y \in C^+\bs\{0\}$ are $\{\R^d, \vn\}$-minimal, not quasi-absorbing and one has $E^\leq(h^+_y) = \{\vn, \R^d, h^+_y\}$.
\end{lemma}

{\sc Proof.} (1) Clearly, closed convex cones and $\vn$ are idempotent. Vice versa, let $\vn \neq A \in \G(\R^d, C)$ be idempotent. On the one hand, if $z \in A$, then $2z \in A$ and hence $nz \in A$ for all $n = 1, 2, \ldots$. Convexity implies $tz \in A$ for all $t > 1$. On the other hand, $z \in A = A \oplus A$ implies that there are sequences $(x_n), (y_n) \subset A$ with $z = \lim (x_n + y_n)$. Convexity of $A$ gives $\frac{1}{2}x_n + \frac{1}{2}y_n \in A$ for all $n$, hence $\frac{1}{2}z = \lim (\frac{1}{2}x_n + \frac{1}{2}y_n) \in A$ by closedness and therefore $tz \in A$ for all $t \geq \frac{1}{2}$ by convexity. Continuing this argument, one gets $tz \in A$ for all $t > 0$.

(2) Let $A, B \in E(\G(\R^d, C))\bs\{\vn\}$, i.e., they are non-empty closed convex cones. Assume $A \supseteq B$. Then $A \oplus B {\f \subseteq} A \oplus A = A$. Since $0 \in B$ one also has $A \subseteq A \oplus B$, thus $A = A \oplus B$ which is $A \leq_H B$ on $E(\G(\R^d, C))$. Conversely, assume $A \oplus B = A$. Since $0 \in A$ one gets $B = \{0\} + B \subseteq A \oplus B = A$.

(3) One can verify $E^\leq(\vn) = \{\vn\}$ and $E^\leq(\R^d) = \{\R^d, \vn\}$ with $\vn \leq_H \R^d$ and hence $\leq_H$ is total on $E^\leq(\R^d)$ (for a definition, compare\cite[Sec. 1.3, p. 13]{Howie95Book}). Moreover, $\R^d$ is absorbing in $(\G(\R^d, C)\bs\{\vn\}, \oplus)$ and $\vn$ is absorbing in $(\G(\R^d, C), \oplus)$.

(4) If $A \in \G(\R^d, C)$ is a non-empty convex cone with $A \supseteq h^+_y$ then either $A = h^+_y$ or $A = \R^d$, hence $h^+_y$ is $\{\R^d, \vn\}$-minimal. Since $y \neq 0$, there is $z \in \R^d$ with $y^\top z = 1$. Then $\{z\} + h^+_y \in \G(\R^d, C)\bs E(\G(\R^d, C))$ and $(\{z\} + h^+_y) \oplus h^+_y = \{z\} + h^+_y\neq h^+_y$, thus $h^+_y \not\in A(\G(\R^d, C))$.
\pend

\begin{corollary}
\label{CorAofWandPofW}
In $(\G(\R^d, C), \oplus)$, one has
\[
A(\G(\R^d, C)) = \{\R^d, \vn\} \quad \text{and} \quad P(A(\G(\R^d, C))) = \{h^+_y \mid y \in C^+\bs\{0\}\}.
\]
\end{corollary}

\begin{theorem}
\label{ThmWandHinG}
For $y \in C^+\bs\{0\}$ and $z \in \R^d\bs\{0\}$ with $y^\top z = 1$ one has
\[
W(h^+_y) = \cb{\{s z\} + h^+_y \mid s \in \R} \cup \cb{\vn, \R^d}.
\]
Moreover, $H(h^+_y) = W(h^+_y)\bs\cb{\vn, \R^d}$ and this set is a group.

Finally, $W(h^+_y)$ is a quasi-minimal ideal in $(\G(\R^d, C), \oplus)$ and $(W(h^+_y), \oplus)$ is a quasi-simple semigroup.
\end{theorem}

{\sc Proof.} It has to be shown that $A \oplus h^+_y = \{s z\} + h^+_y$ for some $s \in \R$ or $A \oplus h^+_y = \vn$ or $A \oplus h^+_y = \R^d$ for each $A \in \G(\R^d, C)$. One has $A \oplus h^+_y = \vn$ if, and only if, $A = \vn$ and $A \oplus h^+_y = \R^d$ if, and only if, $\inf_{a \in A} y^\top a = -\infty$. In the remaining case, $\inf_{a \in A} y^\top a = s \in \R$ holds true and thus $A \oplus h^+_y = \{s z\} + h^+_y$.

Clearly, $W(h^+_y)\bs\cb{\vn, \R^d}$ is a group and therefore, $H(h^+_y) = W(h^+_y)\bs\cb{\vn, \R^d}$ by Theorem \ref{ThmHandAEq}.
\pend

\begin{remark}
\label{RemConlinSpace}
On $\G(\R^d, C)$, a multiplication $\cdot $ with non-negative numbers can be defined as $s \cdot A = \{sa \mid a \in A\}$ (in particular,  $s \cdot \vn = \vn$) for $s > 0$ extended by $0 \cdot A = C$. Using this multiplication along with an obvious extension to negative numbers, one can show that $H(D)$ with $\oplus$ carries a linear space structure whenever $D \in \G(\R^d, C)$ is an idempotent element. In fact, one has
\[
H(D) = \cb{\{z\} + D \mid z \in \R^d}.
\]
The algebraic structure of $(\G(\R^d, C), \oplus, \cdot)$ has been dubbed `conlinear space' in \cite{Hamel05Habil}.
\end{remark}

In set optimization applications also the sets $\F(\R^d, C) = \{A \subseteq \R^d \mid A = \cl(A + C)\}$ with addition $\oplus$ and $\P(\R^d, C) = \{A \subseteq \R^d \mid A = A + C\}$ with element-wise addition $+$ play a role; their semigroup structure is more involved because of the missing convexity.

Finally, one may note that the concepts and arguments in this section straightforwardly extend to a separated, locally convex topological linear space $Z$ with topological dual $Z^*$ instead of $\R^d$.

\section{Perspectives}

First, one may ask if concepts like $A$-minmal ideals and $A$-simple semigroups can be reduced to known cases via the Rees congruence generated by  $A$. In general, the answer is no, but there are important special cases which could be treated in this way. A detailed discussion can be found in \cite{HagerHamelHeyde23ArX} where in particular the importance of bottleneck ideals is emphasized.

Secondly, a more fundamental issue is the following. In the ordered semigroup $(\G(\R^d, C), \oplus, \supseteq)$, a quasi-primitive element $h^+_y$ generates a total order via
\[
A \leq B \quad \Leftrightarrow \quad A \oplus h^+_y \supseteq B \oplus h^+_y.
\]
for $A, B \in \G(\R^d, C)$. Under what conditions to an ordered semigroup does this remain true, i.e., do the quasi-primitive elements generate total orders? Note that in this case the Rees order on $E(\G(\R^d, C))$ conincides with the original order on $\G(\R^d, C)$--except if the absorbing element $\vn$ is involved.

Moreover, for each $A \in \G(\R^d, C)$, one has
\[
A = \sup_{h^+_y \in P(A(\G(\R^d, C)))} A \oplus h^+_y.
\]
The proof of this fact, a basic duality theorem, relies on a separation argument in locally convex spaces. It might be reasonable to ask under which conditions to the (ordered) semigroup and the quasi-absorbing elements such a theorem remains valid in a semigroup set-up.


\end{document}